\documentclass[10pt]{amsart}
\usepackage{graphicx,amssymb,amsmath}
\newtheorem{thm}{Theorem}[section]

\newtheorem{lem}[thm]{Lemma}
\newtheorem{prop}[thm]{Proposition}
\newtheorem{conj}[thm]{Conjecture}
\theoremstyle{definition}

\theoremstyle{remark}
\newtheorem{rem}[thm]{Remark}
\numberwithin{equation}{section}
\usepackage{amsmath,amssymb,amsthm}

\theoremstyle{definition}

\newtheorem{example}[thm]{Example}

\newcommand\Z   {{\bf Z}}

\newcommand\R   {{\bf R}}
\newcommand\C   {{\bf C}}
\renewcommand\P  {{\bf P}}

\newcommand\Pic   {{\rm Pic}}

\newcommand\tr    {{\rm tr}}
\newcommand\Adj   {{\rm Adj\,}}
\newcommand\Aut   {{\rm Aut\,}}
\newcommand\Int      {{\rm K}}

\newcommand\calK      {{\mathcal K}}
\newcommand\calD      {{\mathcal D}}
\newcommand\Kahler    {K\"ahler\ }
\newcommand\CY        {Calabi--Yau\ }

\renewcommand\Im         {{\rm Im}}

\newcommand\im          {{\rm im}}
\renewcommand\Int          {{\rm Int}}

\newcommand{\beas}{\begin{eqnarray*}} 
\newcommand{\eeas}{\end{eqnarray*}} 

\begin{document}

\title{Asymptotic curvature of moduli spaces for Calabi--Yau threefolds}
\author{Thomas Trenner and P.M.H. Wilson} 
\address{Department of Pure Mathematics, University of Cambridge,
16 Wilberforce Road, Cambridge CB3 0WB, UK}
\email {t.trenner@dpmms.cam.ac.uk, pmhw@dpmms.cam.ac.uk}

\date{16 July 2010}

\begin{abstract}
Motivated by the classical statements of Mirror Symmetry, we study certain \Kahler metrics
on the complexified \Kahler cone of a \CY threefold, 
conjecturally corresponding to approximations to 
the Weil--Petersson metric near large complex structure limit for the mirror.  In particular, the naturally defined Riemannian metric (defined via cup-product) on a level set of the \Kahler cone is seen
to be analogous to a slice of the Weil--Petersson metric near large complex structure limit.  
This enables us to give counterexamples to a conjecture of Ooguri and Vafa that 
the Weil--Petersson metric has non-positive scalar curvature in some neighbourhood of the 
large complex structure limit point.\\
\\
\bf Keywords \rm   Mirror symmetry, Weil--Petersson metric, large complex structure limit points, large radius limit points, curvature\\
\\
\bf Mathematics Subject Classification (2000) \rm  Primary 14J32, Secondary 32Q25, 53A15
\end{abstract}


\maketitle

\section*{Introduction}

In this paper, we aim to 
understand the asymptotic behaviour of the Weil--Petersson metric near  large complex structure limit points (defined in terms of maximally unipotent monodromy) on the complex moduli space of \CY threefolds, by using a classical form of mirror symmetry 
and calculating the curvature of certain \Kahler metrics on the complexified \Kahler cone of the mirror. 

This is intimately connected with the theory developed in \cite{Wilson1}.  If $V$ is a \CY threefold with 
$h^{2,0} =0$ and $h^{1,1} =r$, we shall denote the cup-product cubic form on $H^2 (V , \R )$ by $f(y_1 , \ldots , y_r )$.  Let $\calK (V) \subset H^2 (V , \R )$ denote the \Kahler cone, and $\calK _1 \subset 
\calK (V)$ the level set given by $f=1$.  It is a consequence of the Hodge index theorem that the 
restriction of $-\frac{1}{6} (\partial ^2 f /\partial y_i \partial y_j )$ to $\calK _1$ defines a natural Riemannian metric on $\calK _1$.  In \cite{Wilson1}, it was argued in Section 1 that this 
Riemannian manifold reflects the asymptotic Weil-Petersson geometry on the moduli space of 
the mirror near large complex structure limit.  This claim is given a 
more precise justification in this paper; in
particular, see Remark 2.7 below.

The history of expectations concerning the curvature of the Weil--Petersson metric on the moduli space of \CY threefolds is marked by unfulfilled hopes, probably over-influenced by the Weil--Petersson metric on the moduli space of curves, which was known to have negative curvature.  In the \CY threefold case, it was claimed in \cite{Tian} that the holomorphic sectional curvatures were negative, and in \cite{Tod} that all the sectional curvatures were negative.  Both these statements were disproved by the calculations of Candelas et al. \cite{Candelas} for the mirror quintic, where the 1-dimensional 
moduli space was shown to have Weil--Petersson curvature tending to $+ \infty$ as one approaches 
the orbifold point.

The mirror quintic moduli space does however have negative curvature as one approaches the large 
complex structure limit point in moduli.  This fact was generalised in \cite{Wang}, 
where it was shown that, 
whenever the complex moduli space of a \CY threefold is 1-dimensional, the Weil--Petersson metric is asymptotic to the Poincar\'e metric near a large complex structure limit point, and in particular has negative curvature.

There was then a folklore expectation that this asymptotic negativity of the Weil--Petersson curvature 
near large complex structure limit should continue to hold for the complex moduli space having 
arbitrary dimension.    This expectation motivated some of the work carried out in \cite{Wilson1}.  A 
weaker version of this expectation was articulated in Conjecture 3 of \cite{OoguriVafa}, where it was 
conjectured that at least the scalar curvature of the Weil--Petersson metric should be non-positive 
near the points at infinity in moduli (see also the evidence quoted in Example (v) in Section 3 of that paper).
The second author pointed out in Section 2 of \cite{Wilson2} that these expectations of asymptotic 
negativity were likely to be false, with 
conjectural counterexamples provided by the mirrors to the smooth \CY Weierstrass fibrations over the Hirzebruch rational surfaces ${\bf F}_0 = \P^1 \times \P^1$, ${\bf F}_1$ and ${\bf F}_2$.  In Section 3 below, we study in detail the case of the Weierstrass fibration over
${\bf F}_2$; this may also be conveniently described as the desingularization of a general hypersurface of degree 24 in $\P (1,1,2,8,12)$, and is a threefold known in the Physics literature as the STU-model.
 In Theorem 3.3, we see that the mirror to this \CY threefold does indeed provide a counterexample to the conjecture in \cite{OoguriVafa}, since the Weil--Petersson scalar curvature is  unbounded above in any neighbourhood of the large complex structure limit point.  Moreover, we observe that the same 
 thing happens for the mirrors to certain other toric hypersurface \CY threefolds, and in Theorem 3.7 we
 extend this yet further to include a different type of \CY threefold.
 
 In Section 1 of this paper, we review various classical statements of mirror symmetry for \CY threefolds.  Motivated by these results, in Section 2 we study 
 certain \Kahler metrics
on the complexified \Kahler cone of a \CY threefold, 
conjecturally corresponding to approximations to 
the Weil--Petersson metric near large complex structure limit for the mirror.  In Section 3, we apply these general results to certain \CY threefolds with $h^{1,1} = 3$, deducing results on the asymptotic Weil--Petersson geometry of their mirrors.  In particular, we prove the existence of  counterexamples 
to asymptotic non-positivity of the scalar curvature, both of the type predicted  in \cite{Wilson2} 
and of a different type.

The first author is supported on a studentship financed by the EPSRC, 
the Cambridge European Trust and
a Trinity Hall Research Bursary.  
The second author wishes to thank Shinobu Hosono for the benefit of a useful conversation and
email in December 2005.

\section{Asymptotic Mirror Symmetry for \CY threefolds}

A classical formulation of the Mirror Symmetry Conjecture 
involves a \it Mirror Map\rm ,
identifying a neighbourhood of a \it large complex structure limit point \rm 
in the complex moduli space of a \CY threefold with a neighbourhood of a
\it large radius limit point \rm in the K\"ahler moduli space of a mirror \CY threefold, 
under which the 
$B$-model correlation functions on complex moduli are identified 
with the $A$-model (quantum corrected) correlation functions on the K\"ahler 
moduli space of the mirror.

The Mirror Symmetry Conjecture stated in this form is carefully described in 
the book \cite{CK}, in particular Chapters 5, 6, 7 and 8.
Implicit in the statement is the idea of a \it large complex structure limit point\rm ,
expressed in terms of a point of maximally unipotent monodromy (plus an integrality 
condition), or more
precisely  the choice of a simple normal crossing compactification of the 
complex moduli space 
around  such a point.  On the other side of the mirror, we have the idea of a \it large 
radius limit point\rm , which involves a choice of smooth birational model $V$ and 
 a  framing for its \Kahler cone.   A \it framing \rm is a choice of an integral basis $T_1 , \ldots , 
T_r$ of the torsion-free part of $H^2(V , \Z)$ generating a simplicial cone $\sigma$ in 
$H^2 (V , \R )$, whose interior $\Int (\sigma )$ is contained in the K\"ahler cone $\calK (V)$.  
We shall denote by $\im (H^2 (V , \Z ))$ the image of the natural map 
$ H^2 (V , \Z ) \to H^2 (V , \C )$.
Throughout we assume that $h^{2,0} (V) =0$, and so $\calK (V)$ is an open convex cone 
in $H^2 (V , \R )$, and $r$ denotes the Hodge number $h^{1,1}(V)$.

The \it complexified K\"ahler cone \rm is defined by 
$$\calK _\C (V) = \{ \omega \in H^2 (V , \C )\ : \ \Im (\omega ) \in \calK (V) \}/ \im (
H^2 (V , \Z )), $$ the elements of this space usually being written as $B + i \omega$,
where $B$ is often called the \it B-field\rm .  
The \it complexified K\"ahler moduli space \rm is defined 
to be $\calK _\C (V)/ \Aut (V)$.

It is explained on pages 128-9 of \cite{CK} how a framing $\sigma$ gives rise to a complex manifold
$$\calD _\sigma  = \{ \omega \in H^2 (V , \C )\ : \ \Im (\omega ) \in \Int ( \sigma ) \}/ \im (
H^2 (V , \Z )) , $$
and a biholomorphism from $\calD _\sigma$ to $(D^*)^r \subset (\C ^* )^r$ given by  
$$t_1 T_1 + \ldots + t_r T_r \mapsto (q_1 , \ldots , q_r)  = (e^{2
\pi i t_1} , \ldots 
, e^{2 \pi i t_r}),$$
with $D \subset \C$ denoting the unit disc, and $D^* $ the punctured disc.
We may partially compactify $(D^*)^r$  to
$D ^r$, and
 the origin will then be referred to as a \it large radius limit point\rm .   For given class 
 $B + i \omega \in \calD _\sigma$, the limit 
 of $B + i t \omega$ as $t\to \infty$ is this large radius limit point.
 Having chosen the 
framing, we  have uniquely defined local coordinates $q_1 , \ldots , q_r$ on 
$\calD _\sigma$, or equivalently coordinates $t_1 , \dots , t_r$ on its universal cover 
$H^2 (V , \R ) + i \, \Int(  \sigma ) $.  Different framings of the \Kahler cone will give rise to equivalent 
large radius limit points.

The first part of the Mirror Symmetry Conjecture says that to each  
large radius limit point of $V$, defined by a framing on the \Kahler cone, there
is a corresponding large complex structure limit point in the complex moduli 
of the mirror $V^\circ$,  and 
certain uniquely defined local coordinates $q_1 ,
\ldots , q_r$ on an open neighbourhood of the large complex structure limit point, 
so that the Mirror Map identifies these coordinates $q_i$ with those 
defined above on $\calD _\sigma$ for $V$.   For a given 
large complex structure limit point, we use the maximally unipotent monodromy 
condition to produce periods $y_i = \int _{\gamma _i} \Omega$ of the holomorphic 
3-form $\Omega$ ($i=0, 1, \ldots , r = h^{1,2} (V^\circ )$), 
with $\gamma _0 , \gamma _1, \dots , \gamma_r$ being part of a symplectic basis of 
$H _3 (V^\circ , \Z)$, where $y_0$ is holomorphic 
at the limit point and $y_1 , \ldots , y_r$ have logarithmic singularities.  The holomorphic
period $y_0$ is well-defined up 
to a constant multiple, and the $y_i$ are well-defined up to ordering and 
the addition of constant multiples of $y_0$.  We may however 
always normalise the 3-form so that $y_0 = 1$; the integrality conjecture (briefly 
alluded to above) is needed in general to get a uniqueness statement for these periods,
namely that  the $y_i$ for $1 \le i \le r$ are uniquely defined modulo ordering 
and the addition of 
integral multiples of $y_0$.  Assuming that $y_0$ has been normalised to be 1, the 
local coordinates required are $q_i = e^{2\pi i y_i}$; having permuted these coordinates 
appropriately, we should obtain the Mirror Map.

We now introduce the Yukawa couplings on the complex moduli side.  For a given system of 
coordinates $z_1 , \ldots , z_r$, and choice of holomorphic 3-form $\Omega$, we 
 define the Yukawa couplings $Y_{ijk}$ to be 
$$ Y_{ijk} = \int _{V^\circ} \Omega \wedge \nabla _{\partial/\partial z_i} 
\nabla _{\partial/\partial z_j} \nabla _{\partial/\partial z_k} \Omega  ,$$  
or in some references (for instance \cite{CK}) with a minus sign in front, 
where $\nabla$ denotes the 
Gauss--Manin connection.  These couplings depend on the 
local coordinates chosen.

Assume now that we have a symplectic basis 
$A_0, A_1,  \ldots , A_r, B_0 , B_1 , \ldots , B_r$ for $H_3 (V^\circ , \Z )$, and a holomorphic 
3-form $\Omega$ with periods $\zeta _0, \zeta _1 , \ldots , \zeta_ r, \break \xi _0, \xi _1 , 
\ldots , \xi_r$; a suitably general choice of $A_0$ ensures that the ratios $z_i = \zeta _i 
/\zeta _0$, for $1\le i \le r$ form a local holomorphic coordinate system (the $\zeta _0 , 
\ldots \zeta _r$ are called \it homogeneous special coordinates\rm ) \cite{BG, BP}.  
If the corresponding dual basis for the torsion free part of $H^3(V ^\circ , \Z )$ is $\alpha _0 , \alpha _1 , \ldots , \alpha _r , \beta _0 , \beta _1 , \ldots ,\beta_r$, 
the $\zeta _i$ and $\xi _j$ are just the coordinates of the class represented by $\Omega$ with 
respect to the given basis.  If we do not normalize $\Omega$, then 
it is a consequence 
of theory from \cite{BG} that  the $\xi _j$ are holomorphic functions of the 
$\zeta _0 , \ldots , \zeta _r$, and that for some holomorphic function $G$ of the $\zeta _i$, we have 
$\xi _j = \partial G / \partial \zeta _j$ for all $j = 0, 1, \ldots , r$ (see \cite{BP}, pp 237-8).  
If we set $(\Omega , \bar \Omega) = -i \int  \Omega \wedge \bar \Omega  >0$, 
it is checked easily that 
$$  (\Omega , \bar \Omega)  = i \sum _j (\bar \zeta _j {\partial G }/{\partial \zeta _j }
- \zeta _j {\partial \bar G }/{\partial \bar \zeta _j } ) > 0.$$

We shall be interested in the \it Weil--Petersson metric \rm on the complex moduli space; among 
various equivalent definitions, this may be defined as the \Kahler metric with \Kahler 
potential $- \log (\Omega , \bar\Omega )$, where the metric is readily seen to be independent 
of the choice of local holomorphic 3-form $\Omega$.  We wish to reinterpret the above in terms 
of \it affine special coordinates \rm $z_i = \zeta _i/\zeta _0$, for $i = 1, \ldots r$, local holomorphic coordinates on the moduli space.  For details on these special coordinates, and derivations of 
the formulae below, see \cite{Strom}.  With this notation, 
$G (\zeta _0 , \zeta _1 , \ldots , \zeta _r )/\zeta_0 ^2$ is a holomorphic function $F(z_1 , 
\ldots , z_r)$ of these coordinates, and the corresponding Yukawa couplings $Y_{ijk}$ 
are given by ${\partial ^3 F}/{\partial z_i \partial z_j \partial z_k}$;  here $F$ is sometimes called the
\it Gauss--Manin holomorphic prepotential\rm .   In \cite{CK}, this is called the Gauss--Manin potential, 
but we use the term prepotential to coincide with the usage from Physics 
\cite{Strom, HKTY1, HKTY2}, and to
distinguish it from the potential function for the \Kahler metric.  
With respect to these coordinates, one sees that the 
Weil--Petersson \Kahler potential may be written as
 
$$ K =   - \log i \big ( \sum_j  (z_ j - \bar z_j )(\partial _j F + \bar \partial _j \bar F) + 2\bar F - 2 F 
 \big ) ,$$
 where $\partial _j = \partial/\partial z_j$ and $\bar \partial _j = \partial/\partial {\bar z}_j$.

We explained above how, near a large complex structure limit point, we could find a well-defined 
set of periods $y_0, y_1, \ldots , y_r$,  on a neighbourhood 
$(D^*)^r$ of the limit point, where $D \subset \C$ is an open disc, with $y_0$ holomorphic on
$D^r$  
and $y_1 , \ldots , y_r$ having logarithmic singularities.  These gave rise to a holomorphic 
coordinate system $q_1 , \ldots q_r$ on some open neighborhood of the limit point in 
complex moduli; we may however 
consider, within such a neighbourhood, small open sets in the (uncompactified) complex moduli space,
on which therefore $y_0 , y_1 , \ldots , y_r$ may be considered as homogeneous special coordinates.
Normalizing the 3-form $\Omega$ so that $y_0 = 1$, we obtain local  affine special coordinates 
$y_1 , \ldots , y _r$ --- equivalently we may regard these as global coordinates on the 
universal covering of  
$(D^*)^r$.  The corresponding Yukawa couplings $Y_{ijk}$ are seen to be  
globally defined holomorphic functions
on the neighbourhood $(D^*)^r$ of the limit point  in complex moduli, 
although the holomorphic 
prepotential $F(y_1 , \ldots , y_r)$ is well-defined only on the universal cover.  The Yukawa 
couplings we have just defined are called the \it normalized Yukawa couplings \rm (cf. \cite{CK}, 
Definition 5.6.3, where they are taken with a negative sign), or the 
\it$B$-model correlation functions\rm .

The other ingredients in classical mirror symmetry are  (quantum corrected) Yukawa couplings on 
the \Kahler side.  These involve the \Kahler class in the definition of the quantum corrections
via Gromov--Witten invariants; the non-quantum part is given by the coefficients $a_{ijk}$ 
in the cup-product cubic form 
$$f(t_1 , \ldots , t_r) = \sum _{i,j,k} a_{ijk} t_i t_j t_k ,$$ 
with respect to the coordinates defined by the framing (see \cite{CK}, Chapter 7).
These also are called the \it $A$-model correlation functions\rm .  The \it Mirror Symmetry 
Conjecture \rm says that, under the mirror map, the $A$-model and $B$-model correlation 
functions are identified.

The case where this can be described explicitly is that of \CY hypersurfaces 
in 4-dimensional toric varieties $\P (\Delta )$ (where $\Delta$ is a reflexive polytope),
or more precisely the  \CY threefolds obtained by an appropriate resolution of 
singularities.
  According to Batyrev's work \cite{Bat1}, 
the mirror $V^\circ $ of such a \CY hypersurface $V$ 
is obtained by passing to the polar (sometimes called dual) 
reflexive polytope $\Delta ^\circ$. 
We assume for simplicity that all deformations of the complex structure on  $V^{\circ}$  
are realised by deforming the defining polynomial of the 
hypersurface in $\P (\Delta^\circ )$; this last condition 
 is equivalent to the property that, for any codimension 2 face $\Theta ^\circ$ of $\Delta ^\circ$, 
with polar dual face
$\hat \Theta ^\circ$ of $\Delta$, either $\Theta^\circ$ or $\hat \Theta ^\circ$ have 
no interior lattice points.
  Having chosen the desingularisation $V$ 
of the hypersurface in $\P (\Delta )$, we can consider framings $\sigma$ of its 
\Kahler cone, where different 
framings then define equivalent large radius limit points.  However, 
each such $\sigma$ defines also a boundary point of the complex moduli space of $V^\circ$
(with $\sigma$ interpreted as part of the GKZ decomposition, and hence also of the secondary fan).
Conjecture 6.1.4 of \cite{CK} says that this boundary point 
is maximally unipotent (and that different
framings give rise to equivalent boundary points), and the authors remark there
that for threefolds the conjecture follows from an assertion  of Givental \cite{Giv}.    In the 
case where the fan consisting of cones (with vertex at the origin) on the faces of $\Delta ^\circ$ 
may be subdivided at 
integral points on the faces of $\Delta ^\circ$ so as to achieve a regular fan (the polytopes 
of types I and II), this Conjecture was explicitly checked in Section 4 of \cite{HKTY2} --- see also 
(3.38) of \cite{HLY}. 

The generators of the dual cone to
$\sigma$ determine natural toric coordinates to the complex moduli space at this limit, 
these coordinates being  monomials in the coefficients of the defining polynomial.  In 
order however to define the correct mirror map, we need to find local coordinates 
at the limit point defined as  previously via periods, and we therefore  need to calculate periods 
of the holomorphic 3-form.  One way to do this is to observe that they satisfy a certain 
set of generalised hypergeometric differential equations, known as the generalised GKZ system, and 
then to use a variant of the classical method of Frobenius to generate solutions 
from the (unique) holomorphic period (see \cite{CK, HKTY2}, and also (3.38) of \cite{HLY}).   
In the toric hypersurface case, 
one can then define the mirror map uniquely without appealing to the integrality conjecture, 
and the toric mirror map so obtained is conjecturally the same as the mirror 
map one obtains via the integrality conjecture.
 The question of ordering the coordinates correctly is dealt with in the 
toric hypersurface case, as the derivative of the mirror map is  
the \it monomial--divisor map \rm 
$H^{2,1} (V^\circ , \C ) \to H^{1,1} (V , \C )$, 
naturally produced by the toric machinery;
for an explicit description of this, see Section 6.3.2 of \cite{CK}.  
As in the general case, we can define the normalized Yukawa couplings or $B$-model correlation 
functions, and the 
 \it Toric Mirror Symmetry Conjecture \rm says that,  
 under the toric mirror map, the $A$-model and $B$-model correlation 
functions are identified.  It
was this form of the conjecture which was used in \cite{HKTY1} to calculate (conjecturally) the 
Gromov-Witten invariants or instanton numbers of certain three-dimensional 
\CY toric hypersurfaces, thus extending the famous calculations of Candelas et al. \cite{Candelas}
for the quintic. The \it Toric Mirror Symmetry Conjecture \rm
 is expected to hold for all toric hypersurface \CY threefolds using the methods of Givental \cite{Giv}.

Both these classical mirror conjectures may be phrased in terms of holomorphic 
prepotentials.  On the complex
side, one has the Gauss--Manin prepotential $F(y_1 , \ldots , y_r )$ 
as defined above on the universal cover of 
an open neighbourhood of the limit point in complex moduli, 
a holomorphic  function of $y_1 , \ldots , y_r$;  on the \Kahler side, one has the Gromov--Witten 
holomorphic prepotential  defined on $H^2 (V, \R ) + i \, \Int (\sigma) $, as in Chapter 8 of \cite{CK}.
Under the map identifying the complex coordinates $y_i$ and $t_i$, the Gauss--Manin 
prepotential should be identified with the Gromov--Witten prepotential, modulo terms 
which are quadratic, linear or constant in the $t_i$ (see \cite{CK}, Corollary 8.6.3).

In this paper however,  we shall only need asymptotic forms of these conjectures, 
that the mirror map is well-defined and  the 
limit of the $B$-model correlation functions $Y_{ijk}$ as one approaches a 
large complex structure limit point in complex moduli will be the coefficients $a_{ijk}$ of the cubic 
cup-product form on the mirror (with respect to the coordinates determined by the corresponding 
framing).  This may be rephrased as saying that
the $B$-model correlation functions near a large complex structure limit
point in complex moduli will correspond, under the (toric) mirror map, to the topological 
(uncorrected) coupling, plus a correction term which is holomorphic in the coordinates
$q_j = e^{2\pi i t_j}$ and takes value zero at the large radius limit.   
If we consider the cup-product cubic form 
$f$ as a cubic polynomial in the complex coordinates $t_i$ on $H^2 (V, \R ) + i\, \Int ( \sigma )$, then
the asymptotic  form of the conjecture says that under the (toric) mirror map, the Gauss--Manin 
holomorphic prepotential will correspond to $f(t_1 , \ldots ,t_r)/6$, modulo 
terms which are quadratic, linear or constant, and a quantum correction term which 
is holomorphic in the $q_i$ with no constant term, and therefore decays 
exponentially in the $t_i$ near large radius limit.  

The Asymptotic Toric Mirror Symmetry Conjecture was explicitly checked in Section 4 
of  \cite{HKTY2} for toric \CY hypersurfaces corresponding to reflexive polytopes with 
the property that all deformations are of polynomial type and 
which are of Type I or II.
A precise summary of the relevant results is given on page 561 of \cite{HLY}, with the 
results following easily once one has the description of the periods, provided by (3.38) 
of the same paper.  

This Asymptotic Mirror Symmetry property is also a limiting form of 
the Toric Residue Mirror Conjecture
of \cite{BM} (see in particular Section 9), proved in the toric complete intersection case in \cite{Karu}.  
Once we have defined the mirror map,  the limits of the Yukawa 
couplings may be calculated with respect to  natural toric coordinates $z_i$ at the 
relevant boundary point
defined by the choice of framing, where $z_i = z_i (q_1 , \ldots ,q_r)$, since an easy application 
of the Chain Rule and Griffiths  tranversality confirms that in the limit, the 
normalized Yukawa couplings coincide with the Yukawa couplings one obtains with 
respect to the tangent vectors $z_i \partial / \partial z_i$.
An alternative approach to all 
this is via the Nilpotent Orbit Theorem, which yields very accurate approximations to the 
periods near the large complex structure limit point, and hence determines the action of monodromy
on the cohomology (cf. also Theorem 5.1 of \cite{GS}).

\section{AMWP metrics on \Kahler moduli}

Let $V$ now denote a \CY threefold, and $\sigma$ a framing on  its \Kahler cone
$\calK (V)$.  As in the previous section, we have a space $\calD _\sigma$, which 
is biholomorphic to $(D^*)^r \subset (\C ^* )^r$ and has $H^2 (V, \R ) + i \, \Int (\sigma)$ 
as its universal cover.  The framing gives rise to coordinates $t_1 , \ldots , t_r$ on 
$H^2 (V, \R ) + i \, \Int (\sigma )$, where we shall set $t_j = x_j + i y_j$, 
 and coordinates $q_j = e^{2 \pi i t_j}$ ($j = 1, \dots , r$) on 
$\calD _\sigma$.

The asymptotic mirror symmetry property, in the form described in
the previous section, leads us to considering  holomorphic 
prepotentials $F(t_1 , \ldots , t_r)$ on $H^2 (V, \R ) + i \, \Int (\sigma )$, of the form 
$$ F(t_1 , \ldots , t_r) = f(t_1 , \ldots , t_r)/6 + \sum a_{lm} t_l t_m + \sum b_k t_k + c + 
h(q_1 , \ldots , q_r ) , $$
where $f$ denotes the real cubic form given by cup-product, and so  $f(t_1 , \ldots , t_r)$
is a polynomial in the complex variables $t_1 , \ldots , t_r$,
and $h$ is a holomorphic function of the $q_j$ which vanishes at the origin (and 
hence, as a function of the $t_j$, decays exponentially as one approaches the large radius 
limit).  Moreover, we assume that this defines a \Kahler metric on some neighbourhood of the 
large radius limit point in $\calD _\sigma$,  with \Kahler potential locally 
given by 
$$ K(t_1 , \ldots , t_r)  =   - \log i \big ( \sum_j  (t_ j - \bar t_j )(\partial _j F + \bar \partial _j \bar F) + 2\bar F - 2 F   \big ) ,$$
where $\partial _j = \partial /\partial t_j$ and $\bar \partial _j = \partial /\partial \bar t_j$.
We see below that this assumption forces the coefficients $a_{lm}$, $b_k$ of the 
quadratic and linear terms to be real, but there is no such implication concerning the 
constant term $c$ (cf. (4.8) of \cite{HLY}).

If we know that the asymptotic mirror symmetry property holds, then we shall be interested 
in the case where $F$ comes from the holomorphic prepotential defined on the mirror, and 
hence the \Kahler metric defined above corresponds 
under the mirror map to the Weil--Petersson metric on some  
neighbourhood of the corresponding large complex structure limit point of the mirror.  We 
shall use this in the next section to deduce  interesting  curvature properties 
of the Weil--Petersson metric for certain \CY threefolds.  In this section however, we shall 
restrict ourselves to considering the \Kahler moduli space, equipped with a \Kahler metric 
defined as above.

A particular special case of such a metric is when we take $ F = f(t_1 , \ldots , t_r ) /6$, 
with the other terms
being taken to be zero.  An easy calculation (observing 
that we only need to check the case when $f$ is a monomial)
verifies then that 
$$i \bigg (  \sum_j  (t_ j - \bar t_j )(\partial _j F + \bar \partial _j \bar F) + 2\bar F - 2 F 
\bigg )  =  8 f(y_1 , \ldots , y_r )/6,$$ 
and so the metric may also be defined by taking 
a \Kahler potential function  $K_0$ given by $ K _0 (t_1 , \ldots , t_r ) = - \log f(y_1 , \ldots , y_r )$.  
Since $K _0$ is independent of the real coordinates $x_1 , \dots , x_r$,
 the same will be true for the metric 
$(g_{i \bar j})$, which will be given by the 
formula  $$4 g _{i \bar j} = - \partial ^2 (\log f)/ 
\partial y_i \partial y_j = ( \partial f /\partial y_i )( \partial f /\partial y_j )/ f^2 - (\partial ^2 f
/\partial y_i \partial y_j )/f .$$  We see below that this is a metric on the whole complexified 
\Kahler cone $\calK _\C (V)$.  Furthermore, if we 
 define the \it index cone \rm to be the open cone $W \subset 
H^2 (V, \R )$ where $f$ is positive and the Hessian matrix $( \partial ^2 f/\partial y_i 
\partial y_j )$ has index 
$(1, r-1)$, then the \Kahler cone $\calK (V)$ is an open subcone of $W$
by the Hodge index theorem, and we can consider 
the potential function 
$K _0$ to be defined on the complexified index cone $(H^2 (V, \R ) + i W )/\im (H^2 (V, \Z))$.

\begin{lem}  The potential function $K _0$ determines a \Kahler metric 
on the whole complexified index cone, and hence on the open subset $\calK _\C (V)$.  

\begin{proof} Since $K _0$ is independent of the coordinates $x_1 , \ldots , x_r$, and the 
corresponding matrix $g_{i \bar j}$ is always real symmetric, we need only demonstrate positivity.
This is now the statement that we have a Riemannian metric on $W$ defined with respect to the coordinates $y_1 , \ldots , y_r$ by the matrix $- \partial^2 (\log (f))/\partial y_i \partial y_j$, and 
that follows from \cite{Loftin}.  
\end{proof}
\end{lem}

\begin{rem}  We note for future use that 
a more precise result is true --- see Lemma 2.4 of \cite{Tot} and Theorem 1 of \cite{Loftin}.  
Note that  the level set 
$M= W_1$ in $W$, given by $f = 1$, is a submanifold of $W$, and that 
the  restriction of $- \partial ^2 f/\partial y_i \partial  y_j$ to $M$ defines a Riemannian metric, 
sometimes called the \it centro-affine \rm metric on $M$ 
(this is where we use the condition on the index; see \cite{Wilson1, Tot} for further details).  It is shown that the cone 
$W$ equipped with the Hessian metric $ - \partial^2 (\log (f))/\partial y_i \partial y_j$
is isometric (up to a scaling) to the manifold $\R \times M$ equipped with  
 the  product metric of the 
standard metric on $\R$ with the centro-affine metric on the level set $M$.
\end{rem}

For reasons which will become clear soon, we shall refer to the \Kahler metric on 
$\calK _\C (V)$ that we have just defined as the \it Asymptotic Mirror Weil--Petersson \rm metric, or 
more concisely the AMWP metric.

Let us return now to the general case for the holomorphic prepotential $F$, including the 
quadratic, linear and constant terms, and quantum correction term $h$.  The following 
lemma is then a straightforward generalization of the previous calculation

\begin{lem}  For the general case, the potential function $K(t_1 , \ldots , t_r )$ reduces to
$$  - \log \big ( 8 f(y_1, \ldots , y_r )/6  - 4 \sum_{l,m} \Im ( a_{lm}) (x_l x_m + y_l y_m) 
- 4 \sum _k \Im (b_k ) x_k - 4 \, \Im ( c) + H \big  ),$$
where $$H(t_1 , \ldots , t_r) = \sum _j 4 \pi y_j \big (\bar q_j (\partial \bar h / \partial \bar q_j ) - 
q_j (\partial h/\partial q_j)\big ) + 2 (\bar h - h).$$
\end{lem}

\begin{prop}  The matrix of functions $(g_{i\bar j} )$, where 
$g_{i \bar j} = \partial ^2 K/ \partial t_i \partial \bar t_j$,  
is periodic in the real coordinates $x_1 , \ldots , x_r$ if and only if the coefficients 
$a_{lm}$ and $b_k$ are all real.

\begin{proof}  If the coefficients are real, then the $g_{i\bar j}$ are clearly periodic, 
since the same is true for $K$.
For the converse, we need to expand $ g_{i \bar j} = \partial ^2 K/ \partial t_i \partial \bar t_j$ in 
terms of the $x_p$ and $y_q$.  For any given fixed choice of 
values for the coordinates $y_1, \ldots , y_r$,
we would obtain periodic functions in the variables $x_1 , \dots , x_r$.  By inspection
however, 
unless the coefficients $\Im (a_{lm})$ and $\Im (b_k )$ are all zero, we observe that 
any $g_{i \bar j}$ may be expressed as 
 a quotient of two functions, which for large values of the $x_p$ 
are dominated by terms which are polynomial in the $x_p$,  
where the degree of the polynomial in the denominator is 
 more than that in the numerator.   This then would contradict periodicity, 
by considering real coordinates $x_p +n$ for $n$ large, and observing that (for any 
fixed choice of values for $y_1 , \ldots , y_r$) the entries in the matrix $( g_{i \bar j})$ 
become arbitrarily small.
\end{proof}
\end{prop}

We see therefore that we are forced to have \it real \rm coefficients for the quadratic and linear 
terms in the holomorphic prepotential,
and that these therefore do not contribute to the potential function $K$.  Moreover, only
the imaginary part of the constant term in the prepotential contributes to $K$.  We may therefore 
 replace our holomorphic prepotential $F$ by one of the form 
$$ F(t _1 , \dots , t_r) = f(t_1 , \ldots , t_r)/6 + i \, \Im (c) + h(q_1 , \ldots , q_r ), $$
without changing the metric.  We may therefore assume that the \Kahler potential 
is of the form $K(t_1 , \ldots , t_r ) = - \log \big ( f(y_1 , \ldots , y_r ) + a + J(t_1 , \ldots , t_r)
\big ) $, 
for some real constant $a$, and a certain function $J$
decaying exponentially as one approaches the large radius limit.  It is then clear that the 
corresponding matrix of partial derivatives $ g_{i \bar j} = \partial ^2 K/ \partial t_i \partial \bar t_j$
is asymptotic to the matrix defining the AMWP metric as one approaches the large radius limit, 
since each entry of the matrix is a quotient, the dominant terms (for large radius limit) of 
both the numerator and denominator being those occurring for the AMWP metric.  It follows 
in particular that for any general choice of holomorphic prepotential $F$ with real quadratic and 
linear terms, we do obtain a \Kahler metric on some neighbourhood of the large radius limit 
point in $\calD _\sigma$.

\begin{lem}  For each quadruple of indices $i,j,k,l$, the entry $R_{i\bar j k \bar l}$ 
in the curvature tensor for the \Kahler metric determined by a general choice of 
allowable holomorphic prepotential $F$ is asymptotic to the corresponding entry of the 
curvature tensor in the AMWP metric.

\begin{proof} Given the form of the metric $g_{i \bar j} = \partial ^2 K/ \partial t_i \partial \bar t_j$, 
this will follow from the 
formula on page 157 of  \cite{KN}  (valid for any \Kahler metric) that 
$$ R_{i\bar{j}k\bar{l}} = \partial^2 g_{i \bar j}/\partial t_k \partial {\bar t}_l 
 - \sum _{d,e} g^{\bar e d} (\partial g_{i \bar e}/\partial t_k )(\partial g_{\bar j d}/\partial \bar t_l ),$$
 where $(g^{i\bar j})$ is the inverse matrix to $( g_{i \bar j})$, so that $\sum _k g^{i \bar k} g_{j \bar k} 
 = \delta _{ij}$.  The claim follows since the formula gives  a 
 sum of terms, all of which are 
quotients, for which 
the dominant terms (for large radius limit) of 
both the numerator and denominator are those occurring for the AMWP metric.
\end{proof}
\end{lem}

\begin{rem}  In the case when the holomorphic prepotential is 
induced via asymptotic mirror symmetry, the AMWP metric is just what it says, and it is convenient 
to use this name in general.  We observe that this metric is certainly independent of any 
choice of framing.  Indeed, there is a certain amount of evidence from the toric case 
(see Section 4 of \cite{HKTY2} and (4.8) of \cite{HLY}) that, for  holomorphic prepotentials
induced via asymptotic mirror symmetry, 
the linear and constant terms should also be 
independent of the framing.
\end{rem}

\begin{rem}
The AMWP metric is invariant under the involution given by ${t_j \to - \bar t_j}$, which 
implies that the submanifold $\calK (V)$ of $\calK _\C (V)$ is the 
fixed locus of an isometric involution, and hence is a totally geodesic submanifold
(see \cite{Kob}, page 59).  In 
particular, for any tangent plane to $\calK (V)$, the sectional curvature of the \Kahler metric 
on $\calK _\C (V)$ is the same as the sectional curvature of the restricted  
metric on $\calK (V)$, which we noted above is (up to scaling) the 
 product of the standard metric 
on $\R$ with the induced (centro-affine) metric on the level set $\calK _1$ given by $f=1$.  
In the case when 
the holomorphic prepotential is 
induced via asymptotic mirror symmetry, the curvature of the level
set $\calK _1$ reflects the asymptotic behaviour of an appropriate slice of the 
Weil--Petersson metric near the corresponding large complex structure limit point for 
the mirror.  This was a claim made in Section 1 of \cite{Wilson1}, but the justification 
just given is somewhat more convincing than that given in \cite{Wilson1}.
\end{rem}

In the case when asymptotic mirror symmetry applies, we may use it to prove a formula for 
the curvature of AMWP metric, via Strominger's formula for the 
curvature of the Weil--Petersson metric on the complex moduli space of the mirror.  
We denote the AMWP metric, with \Kahler potential $-\log f(y_1 , \ldots , y_r)$,
by the matrix $(g_{i \bar j})$ (with respect to the coordinates $t_j = x_j + i y_j$),  
and let $f_{ijk}$ denote the third partial derivatives of the cubic form $f$. 
On the mirror side, there is a well-known formula, re-proved by Strominger, for the 
curvature of the Weil--Petersson metric \cite{Strom, LuSun, Schu, Wang}.  The 
$(\Omega , \bar \Omega)$ term in the denominator of Strominger's formula 
corresponds asymptotically (as checked above) to $8 f(y_1 , \ldots , y_r)/6$, 
and so the formula predicts an analogous  formula for the curvature 
tensor of the AMWP metric $g_{i \bar j}$.

\begin{conj}  The curvature tensor is given by the formula 
$$ R_{i\bar{j}k\bar{l}} = g_{i\bar{j}}g_{k\bar{l}} +
g_{i\bar{l}}g_{k\bar{j}} -  \sum_{p,q} \frac { g^{p\bar{q}} f_{ikp}f_{jlq}}{64 \, f^2},$$
where $(g^{p \bar q})$ denotes the inverse matrix to $(g_{i \bar j })$.
\end{conj}

This conjecture is an algebraic identity, which should therefore hold
 for arbitrary cubics $f$, on their complexified index cone.  
 We believe that it follows in general from the theory of projective special K\"ahler manifolds 
 \cite{Freed}, but in this paper  we shall only need it for a general ternary cubic,
  i.e. when $r=3$, for which a more elementary proof suffices.
 
 \begin{lem} Suppose $f$ is a ternary cubic, then 
 the above formula for the curvature tensor of the AMWP metric on the 
 complexified index cone holds.

 \begin{proof}  Suppose that $r$ is arbitrary and we make a linear change of variables 
 $y_i = \sum a_{ij} y' _j$, with $A = (a_{ij})$ 
a real invertible matrix; this corresponds to a linear change of variables
in the complex coordinates $t_i = \sum a_{ij} t' _j$.  
Therefore $\partial t_e /\partial t'_j = a_{ej}$ is real for all $e,j$.
We observe now that the components of the 
curvature tensor transform via
$$ R'_{i\bar{j}k\bar{l}} = \sum a_{di} a_{ej} a_{mk} a_{nl} R_{d\bar{e}m\bar{n}}.$$ 
 We claim  that the right-hand side of the desired identity transforms in the 
 same way.  For this we note that 
 $$ g' _{p' q'} = - \frac{\partial }{\partial t'_{p'}} \frac{\partial }{\partial \bar t'_{q'}} 
 \log f(y_1 , y_2 , y_3)  
= \sum \frac{\partial t_p }{\partial t'_{p'}} \frac{\partial \bar t_q}{\partial \bar t'_{q'}} g _{p \bar q}
 = \sum a_{p p'} a_{q q'} g_{p \bar q},$$ 
 or in matrix notation that $g ' = A^t g A\,$, and hence $(g')^{-1} = A^{-1} g^{-1} (A^t)^{-1}$.  It follows 
 therefore that for any holomorphic functions $u$ and $v$,  
 $$ \sum (g')^{p' \bar q'} \frac{\partial u}{\partial t'_{p'}} \frac{\partial \bar v}{\partial \bar t'_{q'}} = 
 \sum g^{pq} \frac{\partial u}{\partial t_p} \frac{\partial \bar v}{\partial \bar t_q}.$$  Applying this 
 to $u = \partial ^2 f /\partial t'_i \partial t'_k$ and $v = 
 \partial ^2  f /\partial  t'_j \partial  t'_l$, where $f = f(t_1 , \ldots , t_r)$ here denotes 
 the cubic considered as a 
 polynomial function in the complex variables, and expressing $u$ and $v$ in terms of the 
 $\partial ^2 f /\partial t_e \partial t_n$, 
 we deduce that the right-hand side of the identity transforms as claimed.  
 
 We now specialise to $r=3$;  
 in order to prove the desired identity for ternary cubics, we may therefore 
assume that the real cubic $f(y_1 , y_2 , y_3)$ is in an appropriate form, 
say Weierstrass canonical form 
$ f(y_1 , y_2 , y_3)  = y_2 ^2 \, y_3  - y_1 ^3  - \lambda  y_1  y_3 ^2  - \mu  y_3 ^3$.  
The validity of the identity in 
this case may be checked by hand (or more conveniently by computer).
 \end{proof}
 \end{lem}
 
If we know  Conjecture 2.8 holds 
in a given case, it has the consequence that the holomorphic 
bisectional curvatures are bounded below by a
fixed constant, namely $-2$; to see this, 
we need only observe that for fixed $i,j$, the real 
symmetric matrix $( f_{ijp}f_{ijq} )$ is positive semi-definite.  This
in turn implies that both the holomorphic sectional curvatures and the Ricci curvatures are 
bounded below (using Section 7.5 of \cite{Zheng}), and hence the same is true for the scalar curvature
(explicitly, the holomorphic sectional curvatures are bounded below by $-2$, the Ricci 
curvatures by $-(r+1)$ and the scalar curvature by $-r(r+1)$); these statements are 
exactly analogous to the ones made concerning the Weil--Petersson metric in  
  \cite{Schu, Wang}.  

\section{Weil--Petersson curvature near large complex structure limits}

If we know that the asymptotic mirror symmetry property holds for a 
pair of \CY threefolds $V^\circ$ and $V$, then the asymptotic behaviour of the Weil--Petersson 
metric near large complex structure limit for $V^\circ$ will be encoded by the behaviour of the 
AMWP metric on the complexified \Kahler cone $\calK _\C (V)$ for $V$.  In Remark 2.7, we observed that for tangent planes to the submanifold $\calK (V)$, the curvature of the AMWP metric is just that of the restricted metric on $\calK (V)$, which (up to a constant scaling) is 
the Hessian metric corresponding to the function 
$- \log f(y_1 , \ldots , y_r )$, where $f$ denotes the cup-product cubic form.  However, in Remark 2.2, 
we commented that this Hessian metric is (up to scaling) 
isometric to the metric on $\R \times \calK _1$ given by the product  of the standard metric on $\R$ 
with the centro-affine metric on 
$\calK _1$ defined by the restriction of the negative Hessian matrix 
$(- \partial ^2 f /\partial y_i \partial y_j )$.  
Modulo a factor of $1/6$, this latter metric is the one studied in \cite{Wilson1}.

We work now under the assumption that the holomorphic bisectional curvatures of the AMWP metric
are bounded below, and hence that the same also 
holds for the holomorphic sectional curvatures.  This holds for instance whenever Conjecture 2.8
is true, which has been shown above to be the case  
 for arbitrary cubics when $r= 3$, and for cubics given by cup-product on \CY threefolds 
when we know that asymptotic mirror symmetry does apply.  We 
remark also that, for any \Kahler metric, if the holomorphic sectional curvatures are bounded 
absolutely, then so too are the sectional curvatures, 
since the sectional curvatures may be expressed in terms of the holomorphic sectional 
curvatures (Lemma 7.19 of \cite{Zheng}).

\begin{prop}  If the holomorphic bisectional curvatures of the AMWP metric on 
$\calK _\C$ are bounded below, and the centro-affine metric on $\calK _1$ has 
sectional curvatures which are unbounded above, then the AMWP metric has Ricci curvatures
 which are unbounded above.  The scalar curvature will then also be unbounded above.

\begin{proof}  From the facts noted above, the holomorphic sectional curvatures of the AMWP metric must also take arbitrarily large values.  Since the Ricci curvature is a sum of holomorphic bisectional 
curvatures (\cite{Zheng}, page 180), the assumptions imply that it too is bounded below by some 
fixed value, but one can find values which are arbitrarily large positive (using the 
unboundedness of the holomorphic sectional curvatures).  Since the scalar curvature is a sum of 
Ricci curvatures, this confirms that the scalar curvature also takes arbitrarily large positive values 
on $\calK _\C (V)$.
\end{proof}
\end{prop}

We now restrict ourselves to the case when $r=3$.  Here, the level set $\calK _1$ is  a 
surface, and a formula was proved in \cite{Wilson1} for the Gaussian curvature $R$
of the centro-affine metric scaled by $1/6$ (which was called the \it Hodge metric \rm on $\calK _1$ 
there), which may be expressed as 
$$ R = -\frac{9}{4} + \frac {S}{4\,  h^2},$$
where $S$ is the $S$-invariant of the cubic $f$, and $h = \det ( \frac{1}{6} \partial ^2 f / \partial y_i \partial y_j)$, 
that is the Hessian determinant  of $f$, scaled by $6^{-3}$.  
The $S$-invariant of a ternary cubic is decribed in \cite{Wilson2, Wilson1, Tot}.
From the formula it follows that if $S>0$ and 
there is a ray in the boundary of the \Kahler cone on which the cubic $f$ does not vanish but the 
Hessian determinant does, then the Gaussian curvature of the centro-affine metric is unbounded above, 
from which it follows as above that the AMWP metric has scalar curvature unbounded above 
on $\calK _\C$.

It was shown in \cite{Wilson} that the boundary of the \Kahler cone is locally polyhedral away 
from the cone given by $f=0$ (known in \cite{Wilson} as the cubic cone).  The codimension one faces 
so obtained correspond to primitive birational 
contractions of the threefold, and  were classified into three types:  
Type I, where only finitely many (rational) curves are contracted; Type II, where a surface (a generalised del Pezzo surface) is contracted to a point; Type III, where a surface $E$ is contracted to a smooth 
rational curve $C$, with $E$ being a conic bundle over $C$.  For rays in the interior of a Type II face, 
it follows immediately that the cubic $f$ does not vanish, but the Hessian 
determinant does.  Here, 
the cubic is of the form $u_1 ^3 + g(u_2 , u_3 )$ with respect to suitable real coordinates 
$u_1 , u_2 , u_3$, from which it follows that the $S$-invariant $S =0$.  
We shall study this case in Theorem 3.7 below.
It should be remarked 
that for all known examples of \CY threefolds with $r=3$, the $S$-invariant of the cubic form is non-negative.

Another case when one has a ray along which the cubic is non-vanishing but the Hessian 
determinant is zero is given by a codimension \it two \rm face corresponding to the contraction of a surface to a point.
This face will be the intersection of faces of Types I or III, and will not automatically imply that $S=0$.
Examples of such \CY threefolds with $h^{1,1} =3$ are given by smooth \CY Weierstrass models over 
smooth surfaces $Y$, where it is shown in \cite{Nakayama} that $Y$ is one of the 
Hirzebruch rational surfaces ${\bf F}_0 = \P ^1 \times \P ^1$, ${\bf F}_1$ or ${\bf F}_2$.  These 
\CY threefolds all 
have Picard groups generated by the class of the section $E \cong Y$ of the fibration and the 
classes pulled back from ${\rm Pic} (Y)$ (generating a hyperplane in $\Pic (V) \otimes \R$).
Here, there will be a codimension two face of the \Kahler cone 
corresponding to the contraction of $E$.  In all three cases, the cubic form is reducible, defining a 
real plane projective cubic consisting of a line (given by the  
classes pulled back from the base) and an irreducible conic intersecting the line in two points;  
the cubic then has $S$-invariant $S>0$, as desired.   For such a \CY threefold, Proposition 3.1 ensures
that the scalar curvature of the AMWP metric is unbounded above;  in Example 3.2 below, we shall 
see this explicitly in the third of the above cases. 
The \Kahler cone is simplicial, and 
the codimension two face of the \Kahler 
cone corresponding to the contraction of the section $E$ is the codimension two face which is 
not contained in the hyperplane of classes pulled back from the base.  In the three cases,  
the ray is the intersection of faces of Types III and III , Types I and III, and Types III and III,  respectively.
In the first case, the faces both correspond to contractions of $E$ (along different rulings), whilst in the 
third case, there are two different exceptional surfaces $E$ and $D$, which intersect along 
the minimal section of $E \cong {\bf F}_2$.

\begin{example}
Let us concentrate on the case above of the \CY threefolds $V$ which are Weierstrass fibrations 
over ${\bf F}_2$. We shall denote the section of the Weierstrass fibration 
by $E$, the pullback of a fibre of the ruling on 
the base by $L$ and the pullback of the $(-2)$-curve on the base by $D$.  These classes 
generate the Picard group of $V$, and the non-zero intersection numbers are given by 
$E^3 = 8$, $E^2\cdot L = -2$, $E\cdot D^2 = -2$ and $E\cdot D\cdot L =1$.  The generators 
of the \Kahler cone are then checked to be
given by $J_1 = E + 2D + 4L$, $J_2 = L$ and $J_3 = D + 2L$.  This leads to the following cubic intersection form:
$$  f = 8y_1^3+12y_1^2 y_3+6 y_1^2 y_2 + 6 y_1 y_3^2 + 6 y_1 y_2 y_3, $$
where the K\"{a}hler cone is given by $y_1, y_2, y_3 \geq 0$.  The $S$-invariant of 
this ternary cubic form takes the value $1$.
In passing, we note that 
the classes $J_1 , J_2 , J_3$ generate $H^2 (V , \Z )$, and so 
there is an obvious natural framing for the \Kahler cone.

The AMWP metric is given by specifying  $- \log f(y_1 , y_2, y_3 )$ to be 
the \Kahler potential.
For the determinant of the metric we then obtain (compare also with Lemma 3.5 below)
$$
    \det(g_{k\bar{l}}) = 
    \frac{27(y_1 y_2+2 y_1 y_3+y_3^2+y_2 y_3)}
    {64 \, y_1^2 (3 y_2 y_3+4 y_1^2+3 y_1 y_2+6y_1 y_3+3 y_3^2)^3}.$$
Using the fact that the Ricci tensor of any \Kahler  metric is given by the formula 
(see \cite{KN}, page 158) 
$$ R_{i \bar j} = - \partial ^2 (\log ( \det(g_{k\bar{l}}))/\partial t_i \partial {\bar t}_j ,$$
we can calculate (using MATHEMATICA) that the scalar curvature is given by 
\begin{eqnarray*}
 & \frac{2}{3 (y_3 (y_2 + y_3) + y_1 (y_2 + 2 y_3))^3} \big ( 
 16 y_1^6 - 9 y_3^3 (y_2 + y_3)^3 + 24 y_1^5 (y_2 + 2 y_3) - \\ &  {}
   27 y_1 y_3^2 (y_2 + y_3)^2 (y_2 + 2 y_3) +
   12 y_1^4 (y_2^2 + 6 y_2 y_3 + 6 y_3^2) - \\ &  {}
   3 y_1^3 (3 y_2^3 + 10 y_2^2 y_3 + 12 y_2 y_3^2 + 8 y_3^3) -
   3 y_1^2 y_3 (9 y_2^3 + 41 y_2^2 y_3 + 64 y_2 y_3^2 +
      32 y_3^3)\big )
\end{eqnarray*}

We can now see explicitly that the scalar curvature blows up to $+ \infty$ if we approach the large 
radius limit along for instance the curve in $\{ {\bf 0} \} \times i \calK (V)$ 
given by $(y_1, y_2, y_3) = (s^2, s, s)$,
letting $ s \to \infty$. It should be noted that using a curve of
the form $(s,a s, b s)$ is not sufficient to get a blowing up of the
scalar curvature, since along such a line it is independent of $s$.

There is however another useful description of this \CY threefold $V$ as the resolution 
of a general hypersurface of degree 24 in weighted projective space 
$\P (1,1,2, 8,12)$; one such hypersurface would be of Fermat type with 
equation $$z_0 ^{24} + z_1 ^{24} + z_2 ^{12} + z_3 ^3 + z_4 ^2 =0.$$ 
The general hypersurface 
inherits singularities from the singular locus of the ambient space, namely it contains 
an elliptic curve  $C$ of $\Z _2$ quotient singularities, which in turn contains an exceptional 
$\Z _4$ quotient singularity (see Section 2 of \cite{HKTY1}).  Resolving singularities, we 
obtain a surface $ D \cong C \times \P^1$ which contracts down to $C$, and a surface 
$E$ isomorphic 
to ${\bf F}_2$ which contracts down to the $\Z _4$ point on $C$; the intersection of these two
surfaces is a curve which is both 
the minimal section of $E$ and a fibre of the ruling on $D$.
 The resolution may be
checked to be a Weierstrass fibration over ${\bf F}_2$, where $E$ and $D$ have the 
same meanings as before.  It is noted in Section 3 of \cite{MorVafa} that such a \CY threefold has
Hodge numbers $h^{1,1} = 3$ and $h^{1,2} = 243$, but that only 242 dimensions of the complex 
moduli are realised as polynomial deformations.

Finally, one should remark that this threefold is often referred to in the Physics literature as the 
STU-model;  when one makes the change of coordinates $U = y_1$, $S = y_2$ and $y_3 = T - U$,
the cubic $f/6$ takes the form $STU + T^2 U + \frac{1}{3} U^3$.
\end{example}

We are now in a position to disprove the conjecture made in \cite{OoguriVafa} concerning the 
non-positivity of the scalar curvature of the Weil--Petersson metric in some neighbourhood 
of the large complex structure limit point.

\begin{thm}  If we consider the mirror $V ^\circ$ to the resolution $V$ 
of a general hypersurface of degree $24$ in 
$\P (1,1,2, 8,12)$, with the natural framing on $\calK (V)$, then there is a
corresponding large complex structure limit point in the moduli space of $V^\circ$, 
and in any neighbourhood of this boundary point the scalar curvature of the 
Weil--Petersson metric is unbounded above.

\begin{proof}  Modulo a few mechanical checks, this follows by combining the various results 
we have proved above.  The weighted projective space $\P (1,1,2, 8,12) = \P (\Delta) $, for 
$\Delta$ the polytope with vertices.
\begin{align*} &v_1=(1,-1,-1,-1) \\
&v_2=(-1,2,-1,-1) \\ &v_3=(-1,-1,11,-1) \\ &v_4=(-1,-1,-1,23) \\
&v_5=(-1,-1,-1,-1). \end{align*}  For the polar polytope $\Delta ^\circ$,  we have 
vertices: \begin{align*} &v_1^* = (1,0,0,0) \\ &v_2^* = (0,1,0,0) \\ &v_3^* = (0,0,1,0)\\ &v_4^* = (0,0,0,1)\\ &v_5^* = (-12,-8,-2,-1).\end{align*} The other integral 
points (apart from the origin) 
of $\Delta^\circ$ are \begin{align*} &v_6^*=(-3,-2,0,0)=\frac{1}{2} v_3^* +
\frac{1}{4} v_4^* + \frac{1}{4} v_5^* \\
&v_7^*=(-6,-4,-1,0)=\frac{1}{2} v_4^*+\frac{1}{2} v_5^* \\
&v_8 ^* = (-1,-1,0,0) = \frac{1}{2} v_1 ^* + \frac{1}{4} v_3 ^* + \frac{1}{8} v_4 ^* + \frac{1}{8} v_5 ^*\\
&v_9^*=(-2,-1,0,0) = \frac{1}{3} v_2 ^* + \frac{1}{3} v_3 ^* + \frac{1}{6} v_4 ^* + \frac{1}{6} v_5 ^*\\
&v_{10} ^* = (-1,0,0,0) = \frac{1}{2} v_9 ^* + \frac{1}{2}v_2 ^*.
\end{align*} \\
Note that we have one interior point of a codimension 2 face of $\Delta ^\circ$, namely 
$v_6^*$. The corresponding dual face of $\Delta$ is spanned by $v_1$
and $v_2$, and there are no integral points in the interior of
this face.  Thus (unlike the case of  $V$), the complex deformations of $V^\circ$ all arise from
deformations of the defining polynomial.  It is readily checked that the fan given by cones 
on the faces of $\Delta ^\circ$ with vertices at the origin may be subdivided at the extra integral
points $v_6^* , v_7^*, v_8^*, v_9 ^*, v_{10}^*$ so as to achieve a regular fan.  
Because of the existence of interior points $v_8 ^*,  v_9 ^*, v_{10}^*$ 
 of  codimension one faces, in the terminology of \cite{HKTY2, HLY} 
the polytope $\Delta ^\circ$ is of Type II rather than Type I.

These calculations ensure that there is a large complex structure limit point for $V^\circ$ 
corresponding to the natural framing on the \Kahler cone of $V$, and that the  
asymptotic toric mirror symmetry property does hold.  Thus the asymptotic form 
the Weil--Petersson metric and its scalar curvature near large 
complex structure limit for $V^\circ$ 
will be encoded by the AMWP metric and its scalar curvature on $\calK _\C (V)$;  as we 
saw in  the previous calculations the
 scalar curvature is unbounded above.
\end{proof}
\end{thm}
 
\begin{rem}  Suppose now we choose a framing $\sigma$ of the above \Kahler cone 
whose interior is contained in the open subcone of $\calK (V)$ where the AMWP 
metric has positive
scalar curvature.  Such a framing will correspond to a 
large complex structure point of the 
mirror (on a different compactification of the moduli space) with the property that, in some
open neighbourhood of this boundary point, the Weil--Petersson metric 
has everywhere positive scalar curvature.
\end{rem}

 By studying the cubic intersection forms for toric hypersurface \CY 
threefolds listed in Appendix C of 
 \cite{HLY}, we find other candidates with $h^{1,1} = 3 $
 where the above theory holds in exactly the same way.
 The reader should be aware that the way these forms are tabulated in the Physics literature 
 omits binomial coefficients, and so what appears as $8 J_1^3 + 2 J_1^2J_2 + 4 J_1^2 J_3 + 
 J_1 J_2 J_3 + 2 J_1 J_3^2$ in Appendix C of \cite{HLY} is exactly the cubic form for the 
 STU-model calculated above.  For all the examples listed there, the \Kahler cone is 
 simplicial with generators $J_1$, $J_2$ and $J_3$.
 For those examples which are of Types I or II, we can 
 check which intersection forms have strictly positive $S$-invariant and which models have 
 codimension two faces (i.e. rays) in the boundary of the \Kahler cone along which the 
 Hessian determinant vanishes but the cubic form is non-zero.  Apart from 
 $V_{24} \subset \P (1,1,2,8,12)$, we find further examples  
 $V_{10} \subset \P (1,1,2,2,4)$, $V_{16} \subset \P (1,1,3,3,8)$,
 $V_{18} \subset \P (1,2,2,4,9)$ and 
 $V_9 \subset \P (1,1,1,2,4)$, where a precisely
 analogous statement to Theorem 3.3 will hold.  As remarked however in Section 5 of \cite{HLY}, 
 the last two of these are in fact isomorphic.

When the $S$-invariant is zero, we do not have any blowing up of sectional curvatures
of the AMWP metric  on
tangent planes to the slice $i \calK (V)$ in $ \calK _\C (V)$; the scalar curvature may still 
however be unbounded above.  To see why this is, let us continue to restrict ourselves to the case 
of $r=3$, and let 
$$ g_{i\bar j} = \big ( (\partial f /\partial y_i ) (\partial f /\partial y_j ) -  
f \, \partial^2 f /\partial y_i \partial y_j \big)/4f^2 $$
denote the AMWP metric.

\begin{lem}  If $G$ denotes the matrix with entries $$G_{ij} = 
(\partial f /\partial y_i ) (\partial f /\partial y_j ) -  
f\,  \partial^2 f /\partial y_i \partial y_j = - f^2 \partial ^2 (\log f )/\partial y_i \partial y_j ,$$ 
then $\det G = \frac{1}{2} f^3 H$.

\begin{proof}  Under linear changes of variables $y_i = \sum a_{ij} y' _j$, with $A = (a_{ij})$ 
a real invertible matrix, we note that the matrix $G$ transforms to $A^t G A$ and 
the Hessian matrix of $f$ transforms similarly.  Hence both $\det G$ and the 
Hessian determinant $H$ transform 
via multiplication by $(\det A)^2$.  Therefore, in order to prove the desired identity, we may 
assume that $f$ is in an appropriate canonical form, say Weierstrass  canonical form 
$ f =  y_2 ^2 \, y_3  - y_1 ^3  - \lambda  y_1  y_3 ^2  - \mu  y_3 ^3$.  The validity of the identity in 
this case may be checked by hand (or computer).
\end{proof} 
\end{lem}

Let us now set $B = \Adj (G)$, the matrix of cofactors of $G$, and let $C(i,j)$ denote the 
rank 1 positive semi-definite matrix with $C(i,j)_{pq} = f_{ijp} f_{ijq}$.  From the formula 
proved in Lemma 2.9 for the curvature of the AMWP metric, we may deduce the following 
criterion for the unboundedness of the scalar curvature.

\begin{lem}
If there exists a ray in the boundary of the \Kahler cone along which the cubic $f$ does not vanish,
the Hessian determinant $H$ does vanish, and for which $\tr (B\, C(i,j))$ does not vanish for some $i,j$, 
then the scalar curvature of the AMWP metric is unbounded above as one approaches the given 
ray.

\begin{proof}
Using Lemma 2.9, we consider the formula for $R_{i \bar i j \bar j}$; 
using the above lemma, the final term in the formula reduces to 
$- \frac{1}{8}  \tr (B\, C(i,j))/H f^3$.  Our assumptions therefore imply that $- R_{i \bar i j \bar j}$ is 
unbounded above, and hence the same is true for the corresponding holomorphic 
bisectional curvature.  Since the holomorphic bisectional curvatures are bounded below, 
it follows immediately that some Ricci curvature is unbounded above, and hence the same is true 
for the scalar curvature.
\end{proof}
\end{lem}

We now return as promised to the case where some codimension one  face of the \Kahler cone corresponding 
to a Type II contraction.

\begin{thm}  If $V$ is a \CY threefold with $h^{1,1} = 3$, for which there is a face of the 
\Kahler cone determining a contraction of Type II, then the scalar curvature of the AMWP 
metric on $\calK _\C (V)$ is unbounded above.

\begin{proof}  By choosing appropriate linear coordinates, we may assume that the cubic $f$ 
takes the form  $f = y_1 ^3 + y_2 ( a y_2 ^2 + 3 b y_2 y_3 + 3 c y_3 ^2)$, 
with $a, b, c$ real coefficients, not both of $b$ and $c$ being zero, and where 
the hyperplane defining the Type II face is $y_1 =0$. 
In the proof of Lemma 2.9, we saw that 
the formula for the entries of the curvature tensor remained valid whatever linear coordinates 
were used.  We therefore use the above coordinates.  With the notation as in Lemma 3.5, we 
set $C = C(1,1)$, and calculate that 
$ \tr (B C)$ is a non-zero constant multiple of 
 $$   ( (b^2 - a c) y_2 ^2 + bc\, y_2 y_3 + c^2 y_3 ^2 )(f - 3 y_1 ^3)f,
$$
 which is non-vanishing at the general point of the hyperplane $y_1 =0$.  
In particular, it follows from Lemma 3.6 that the scalar curvature of the 
AMWP metric is unbounded above.
\end{proof}
\end{thm}

\begin{example}  Consider a \CY threefold $V$ from \cite{HLY} given as a resolution of a hypersurface 
$V_{16} \subset \P (1,1,1,5,8)$.  It is shown in \cite{HLY} that this is of Type II, and so 
the asymptotic toric mirror symmetry property holds, and $ h^{1,1} =3$.  The topological cubic form is  
 $$f = 50 y_1^3 + 30 y_1 ^2 y_2 + 6 y_1 y_2^2 + 240 y_1 ^2 y_3 + 96 y_1 y_2 y_3 + 
9 y_2 ^2 y_3 + 384 y_1 y_3^2 + 75 y_2 y_3 ^2 +203 y_3^3 .$$
The Hessian of this vanishes along the hyperplane $y_2 = 0$, and this hyperplane defines 
a contraction of $V$ of Type II.  In particular, the $S$-invariant is zero.
Theorem 3.7 shows that the scalar curvature of the AMWP metric on $\calK _\C (V)$ is unbounded 
above, and hence, by asymptotic mirror symmetry, the scalar curvature of the Weil--Petersson metric is 
unbounded above in any neighbourhood of the corresponding large complex structure limit point of 
the complex moduli space of the mirror.

Similar statements are true for the \CY threefold arising from $V_{12} \subset \P (1,1,1,3,6)$, which is also of Type II (and so the asymptotic toric mirror symmetry property holds)
 and has $ h^{1,1} =3$.  The topological cubic form is
$$ f = 18 y_1 ^3 + 18 y_1 ^2 y_2 + 54 y_1 ^2 y_3 + 6 y_1 y_2 ^2 + 36 y_1 y_2 y_3 + 54 y_1 y_3^2 
+ 3 y_2^2 y_3 + 9y_2 y_3^2 + 9 y_3^2,$$
and the Hessian vanishes along the codimension one face of the \Kahler cone  given by $y_2 =0$.
\end{example}

In a similar way, certain toric hypersurface Calabi--Yau threefolds admitting a Type II contraction 
and with $h^{1,1} =2$ also provide counterexamples to the original conjecture.

\bibliographystyle{ams}

\begin{thebibliography}{20}
\bibitem {Bat1}Batyrev, V.V.: Dual polyhedra and mirror symmetry
for Calabi-Yau hypersurfaces in toric varieties.  J. Alg. Geom. {\bf 3},
493-535 (1994)
\bibitem {BM}Batyrev, V.V., Materov, E.N.: Toric Residues and Mirror Symmetry,  Mosc. Math. J. 
\bf 2\rm , 435-475 (2002)
\bibitem {BP}Bertin, J., Peters, Ch.: Variations de structures de Hodge, vari\'et\'es de Calabi--Yau 
et sym\'etrie miroir.  In:  Introduction \`a la th\'eorie de  Hodge (J. Bertin, J.-P. Demailly, 
L. Illusie, Ch. Peters), pp 169-256.  Soci\'et\'e Math\'ematique de France, 1996
\bibitem {BG}Bryant, R., Griffiths, P.: Some observations on the infinitesimal period relations for regular threefolds with trivial canonical bundle, Arithmetic and Geometry vol. II, M. Artin and J. Tate (eds.), Progress in Math. vol. 36, Birkh\"{a}user, Boston, 1983, 77-ö102.
\bibitem {Candelas}Candelas, P., De la Ossa, X. C., Green, P. S., Parkes, L.: An exactly soluble superconformal theory from a mirror pair of Calabi-Yau manifolds, Phys. Lett. B.
{\bf 258},  118--126 (1991)
\bibitem {CK}Cox, D.,  Katz, S.: Mirror
Symmetry and Algebraic Geometry. 
Mathematical Surveys and Monographs, {\bf 68}.  AMS Publishing, 1999
\bibitem{Freed}  Freed, D.S.:  Special K\"ahler manifolds, Comm. Math. Phys. {\bf 203},  31-52 (1999)
\bibitem{Giv}  Givental, A.:  A mirror theorem for toric complete intersections. In: Topological field theory, primitive forms and related topics (eds Kashiwara, M. et al.).  Proceedings of the 38th Taniguchi symposium, Kyoto, Japan, December 9-13, 1996.  Birkh\"auser Prog. Math. \bf 160\rm , 141-175 (1998)
\bibitem {GS}  Gross, Mark; Siebert, Bernd:  Mirror Symmetry via logarithmic degeneration data II.  
arXiv:07092290
\bibitem {HKTY1}Hosono, S., Klemm A., Theisen S., Yau, S.-T.: Mirror Symmetry,
Mirror Map and Applications to Calabi-Yau Hypersurfaces, Comm. Math.
Phys. \bf 167\rm , 301-350 (1995)
\bibitem {HKTY2}Hosono, S., Klemm A., Theisen S., Yau, S.-T.: Mirror Symmetry,
Mirror Map and Applications to Complete Intersection Calabi-Yau Spaces, Nuclear Physics B. 
\bf 433\rm ,  501-552 (1995)
\bibitem{HLY}Hosono, S., Lian, B.H., Yau, S.T.:   GKZ-Generalized Hypergeometric Systems 
 in Mirror Symmetry of \CY Hypersurfaces.  Commun. Math. Phys. \bf 182\rm , 535-577  (1996)
 \bibitem {Karu} Karu, K.:  Toric residue mirror conjecture for \CY complete intersections, J. Algebraic Geom. \bf 14\rm , 741-760 (2005)
 \bibitem{Kob}Kobayashi, S.:  Transformation Groups in Differential Geometry.  
 Berlin-Heidelberg: Springer 1985 
\bibitem {KN}Kobayashi, S., Nomizu, K.: Foundations of Differential
Geometry, Volume II.  John Wiley and Sons, 1969
\bibitem {Loftin}Loftin, J.C.: Affine Spheres and K\"{a}hler-Einstein Metrics, Math. Res. Letters,  
\bf 9\rm , 425--432 (2002)
\bibitem {LuSun}Lu, Z., Sun, X.: Weil-Petersson Geometry on the
Moduli Space of Polarized Calabi-Yau Manifolds, J. Inst. Math.
Jussieu, \bf 3\rm , 185--229 (2004)
\bibitem {MorVafa}Morrison, D., Vafa, C.: Compactifications of F-Theory on Calabi--Yau Threefolds -- I, Nuclear Physics B. \bf 473\rm ,  74--92 (1996)
\bibitem{Nakayama} Nakayama, N.: On Weierstrass models.  In: Algebraic geometry and commutative algebra, Vol.II, pp 405-431.  Kinokuniya, Tokyo, 1988
\bibitem {OoguriVafa} Ooguri, H., Vafa, C.: On the geometry of the string landscape and the swampland, Nuclear Physics B. \bf 766\rm ,  21--33 (2007) 
\bibitem{Schu}  Schumacher, G.:  The curvature of the Petersson--Weil metric on the moduli space of 
K\"ahler--Einstein manifolds.  In: Complex Analysis and Geometry 
(eds V. Ancona \& A. Silva), pp 339-354.  Plenum Press, New York 1993
\bibitem{Strom}Strominger, A.: Special Geometry, Comm. Math.
Phys. \bf 133\rm , 163--180 (1990)
\bibitem {Tian}Tian, G.: Smoothness of the Universal Deformation
Space of Compact Calabi-Yau Manifolds and its Petersson-Weil metric.
In S.-T. Yau, ed., Mathematical aspects of string theory, vol. 1,
629--646. World Scientific, 1987
\bibitem {Tod}Todorov, A.N.: The Weil-Petersson geometry of the 
moduli space of $SU(n\leq 3)$ (Calabi-Yau) manifolds. I. Comm. Math. Phys \bf 126\rm ,
 325--346 (1989)
\bibitem {Tot}Totaro, B.:  The Curvature of a Hessian Metric.  Int. J. Math. \bf 15\rm , 369-391 (2004)
\bibitem {Wang}Wang, C.-L: Curvature Properties of The Calabi-Yau
Moduli, Documenta Mathematica \bf 8\rm , 577--590 (2003)
\bibitem{Wilson}  Wilson, P.M.H.:  The \Kahler cone on \CY threefolds, Invent. math. \bf 107 \rm (1992), 561-583; Erratum: Invent. math. \bf 114\rm , 231-232 (1994)
\bibitem {Wilson1}Wilson, P.M.H.: Sectional Curvatures of K\"{a}hler
moduli, Math. Ann. \bf 330\rm , 631--664 (2004)
\bibitem {Wilson2}Wilson, P.M.H.: Some Geometry and Combinatorics for the $S$-Invariant of Ternary Cubics, Experiment. Math. \bf 15\rm , 479--490 (2006)
\bibitem{Zheng}Zheng, F.:  Complex Differential Geometry.  Studies in Advanced Mathematics {\bf 18}.
AMS Publishing and International Press, 2000.

\end{thebibliography}

\end{document}